# Constructivism versus Cognitive Load Theory: In Search for an Effective Mathematics Teaching

2 authors, including:

Bustang Buhari
Universitas Negeri Makassar
4 PUBLICATIONS   31 CITATIONS

SEE PROFILE

Some of the authors of this publication are also working on these related projects:

Mathematics modelling as Instruction Approaches   View project

Investigating Indonesian secondary school students' and teachers' probabilistic reasoning   View project

# Constructivism versus Cognitive Load Theory: In Search for an Effective Mathematics Teaching

Hamzah Upu[1]
Bustang

Department of Mathematics
State University of Makassar
Jl AP Pettarani, Makassar, 90224, Indonesia

**Abstract**
Two major learning theories have dominated recent literature on optimizing knowledge acquisition: constructivism and cognitive load theory. Constructivism, on the one hand, gives preeminent value to the development of students' self-regulated process of constructing mathematical concepts. Its basic tenet is that students acquire their own mathematical understanding by constructing them from the inside rather than by internalizing them from the outside. Cognitive load theory, on the other hand, suggests that the free exploration of a highly complex environment may cause a heavy working memory load and led to poorer learning. Advocates of this view further argue that constructivist strategies provide learners with information that exceeds their working memory capacity, and thus fail to efficiently guide learners' acquisition of mathematical knowledge. The current study describes the elements of constructivism theory and their cognitive basis and show how they can be aligned with the structures that constitute human cognitive architecture. More specifically, we present several ways in which cognitive load can be managed by these elements and so facilitate mathematical learning.

*Keywords:* constructivism, cognitive load theory, mathematical learning

**Introduction**
Learning is much more than memorizing. Learning refers to the acquisition of knowledge through interactions with, and observation of, the physical word and the creatures that inhabit it (Ashman & Conway, 1997). In order to really understand and be able to apply knowledge, students must work to solve problems, to discover things for themselves, and to struggle with ideas. The question of how to help students learn particular knowledge, skills, and concepts that will be useful in their life is at the core of the argument presented by Kirschner, Sweller, and Clark (2006). The authors compare minimally guided instructions with instructional approaches that provide direct instructional guidance of the student learning process. They define minimally guided instruction as "one in which learners, rather than being presented with essential information, must discover or construct essential information for themselves" and then inversely define direct instruction as "providing information that fully explains the concepts and procedures that students are required to learn as well as learning strategy support that is compatible with human cognitive architecture" (p. 1).

---

[1] Correspondence should be addressed to Hamzah Upu, State University of Makassar, Department of Mathematics, Makassar, Indonesia. E-mail: hamzahupu@gmail.com

In their argument, Kirschner, Sweller, and Clark (2006) affirm that minimal guided instruction approaches are less effective and efficient than fully guided instruction approaches because they ignore the structures that constitute human cognitive architecture. On the contrary to this, they put a strong emphasis on direct, strong instructional guidance, as an effective and efficient way to teach students. By referring to several studies concerning the efficacy of direct instruction (e.g., Klahr & Nigam, 2004), they claim that students learn more deeply from strongly guided instruction than from constructivist or discovery approaches. Opposing this claim, Kuhn and Dean (2006) have found that direct instruction does not work so well for robust acquisition or for maintenance knowledge over time.

According to Hmelo-Silver et al. (2007), there are two major flaws of Kirschner and his colleague's argument. The first is in their pedagogical point of view. Kirschner et al. (2006) have included several distinct pedagogical approaches – constructivist, discovery, problem-solved, experiential, and inquiry-based teaching – under their "minimally guided" umbrella. In agreement with to Hmelo-Silver et al (2007), I argue that some of these approaches, in particular constructivist, cannot be equated with minimally guided instruction. In contrast to Kirschner et al's point of view, I assert that the elements of constructivist approaches allow for flexible adaptation of guidance, making these instructional approaches more compatible with the manner in which our cognitive structures are organized. The second flaw in Kirschner et al. as identified by Hmelo-Silver et al. is in their evidentiary base. The claim of Kirschner and his colleagues that constructivist approaches are ineffective contrasts with to empirical evidence that indeed support the efficacy of constructivist as instructional approaches (e.g., Sultan, Woods, & Koo, 2011; Tatli & Ayas, 2012; Blink, 2000). This evidence suggests that constructivist approaches can foster deep and meaningful learning as well as critical thinking of the students.

In the work presented here, I will discuss how constructivist approaches may provide instructional guidance and evidence that encourage the effectiveness of these pedagogical approaches and cannot be categorized as a minimally guided instruction. Furthermore, I will also describe the elements of constructivist approaches and show how they align to the structures that constitute human cognitive architecture.

**Constructivist Approaches are not Minimally Guided Instructions**

Constructivism has been a leading if not the dominant theory or philosophy of learning, since it has been applied in almost all subjects (Noddings, 1999). This theory implies a new kind of pedagogy where the emphasis will be more on what students do than what teachers do, and where there will be performance assessment of student learning rather than standardized achievement testing (Iran-Nejad, 2001). The essence of constructivist theory is the idea that knowledge is not transmitted directly from one knower to another or from teacher to the students. Rather, knowledge is actively built up by the learners in their minds. As one of the leading exponents of constructivism, von Glasersfeld (2000) said:

> It holds that knowledge is under all circumstances constructed by individual thinkers as an adaptation to their subjective experience. This is its working hypothesis and from it follows that for a constructivist there cannot be anything like a dogmatic body of unquestionable knowledge. The task is to show that and how what is called knowledge can be built up by individual knowers within the sensory and conceptual domain of individual experience and without reference to ontology. (p. 4)

This view has profound implications for teaching, as it suggests that students, instead of the teacher, organize information, explore the learning environment, conduct learning activities, and understand their own learning. Furthermore, constructivist approaches do not relinquish teacher to control of the classroom, as is typically implemented, neither does a sense of ownership mean collecting students' ideas and teaching them back to them (Iran-Nejad, 2001). Instead, the teacher in the constructivist learning environment, facilitates the construction of knowledge of students by teaching in ways that make information meaningful and relevant to students, by giving them opportunities to discover or apply ideas themselves, and by teaching students to be aware of and consciously use their own strategies for learning.

The constructivist approaches work in exactly in opposite order as the traditional approach by starting with problems, then the teacher help them figure out how to do the operations. Moreover, students in constructivist learning environment, instead of construct their knowledge by themselves, are given scaffolds and have to climb these scaffolds by themselves in order to reach higher understanding.

I agree with Vogel-Walcutt et al. (2011) that several instruction approaches fall under the umbrella of constructivism such as problem-based learning (Schmidt et al., 2007), discovery learning (Dean & Kuhn, 2006) and inquiry learning (Hmelo-Silver, Duncan & Chinn, 2007). This implies that problem-based learning and inquiry learning are categorized as an example of constructivist approaches.

**The Use of Scaffolding in Constructivist Approaches**
As I have examined the broad variety of constructivist approaches, there are several key elements that make these approaches cannot be categorized as a minimally guided instruction. The claim by Kirschner et al. (2006) that constructivist approaches provided minimal information to the learners are contradicted with the point of view of advocates of constructivist approaches (e.g., Blikn, 2000; Simon, 1995).

In a constructivist's point of view, the students are not expected to construct everything on their own. In constructivist approaches such as problem-based learning and inquiry learning, students are provided with so-called scaffolding to support students' learning of both how to do the task as well as why the task should be done that way (Hmelo-Silver, 2007). The term scaffolding is defined by Reiser (2004) as the process by which the teacher or more knowledgeable peer assist a learner, altering the learning task so the learner can solve problems or accomplish tasks that would otherwise be out of reach. Quintana et al. (2004) conceived scaffolding as a key element of cognitive apprenticeship, whereby students become increasingly accomplished problem-solvers given structure and guidance from mentors who scaffold students through coaching, task structuring, and hints, without explicitly giving students the final answers. In their research, Quintana et al. (2004) design three constituent processes of scaffolding, that is

> Sense making, which involves the basic operations of testing hypothesis and interpreting data; process management, which involves the strategic decision involved in controlling the inquiry process; and articulation and reflection, which is the process of constructing, evaluating, and articulating what has been learned. (p. 341)

According to several studies implementing scaffolding to support students' learning (e.g., Puntambekar & Kolodner, 2005; Quintana, 2004), there is evidence that scaffolding makes learning more tractable for students and helps them to deal with complex problems. Moreover, scaffolding is not only focused on interaction with teacher or peer, as the source of assistance, but also with the technology design in which technological tools provide some types of assistance. One of example in this manner is the research done by Kim and Hannafin (2011), show that technology-enhanced scaffolds are effective in supporting scientific inquiry learning.

**Human Cognitive Architecture**

A brief explanation of human cognitive architecture and the fundamental elements of cognitive load theory is presented in this text in order to give a picture for the readers about those two terms since they are at the core of the Kirschner et al.'s argument. According to Kirschner et al. (2006), human cognitive architecture determines the manner in which our cognitive structures are organized. Most cognitive theories treat human cognitive architecture by using long-term and short-term or working memory (Kirschner, 2002).

Short-term or working memory is defined as the memory that is used for all conscious activities such as reading the text and it is the only memory that we can monitor (Sweller, 2004; Kirschner, 2002). The problem, especially in relating to instructional design, is that the capacity of working memory in saving information is limited. Miller (1994) in his research has found that working memory only can hold about five to nine items or elements of information that have not been previously learned or known, in a certain time.

Long-term memory, in contrast, is defined as the repository that consists of large and relatively permanent store of information, knowledge and skills (Sweller, 2004; Kirschner, 2002). Most of the cognitive scientists believe that long-term memory can hold unlimited amounts of information including large, complex interactions and procedures. In relation to this, Kirschner (2002) and Sweller (2004) have argued that instructional designers have to consider how is the information stored and organized in long-term memory so that the learners can access this information whenever they need it.

**Cognitive Load Theory**

Cognitive load theory is based on the assumption that human cognitive architecture can only process a limited amount of information in working memory in a certain time (Kirschner, 2002). As a result, any information, knowledge or skills that is presented to the learners and exceeds this capacity may only enter working memory, but will not be stored into long-term memory. Cognitive load theory is concerned with the limitation of working memory capacity and the manner in which the level of cognitive load can be measured to promote an effective learning. Thus, the proponents of cognitive load theory (e.g., Kirschner et al., 2006) state that the aim of all instruction is to alter long-term memory. The best way to achieve this goal is by using direct instructional guidance in which the learners are strongly guided by providing information that fully explains the concepts and procedures that are required to learn.

Three discrete types of cognitive load have been defined (Kirschner, 2002; Sweller, 2004) namely intrinsic, extraneous, and germane cognitive load. The intrinsic cognitive load is affected by the learning content of the task of subject matter itself. Intrinsic cognitive load takes place in the mind of learners when the elements of the to-be-learned material are highly interconnected (Kirschner, 2002). The extraneous cognitive load deals with the manner in which the task information is presented to learners and also the learning activities required of learners. More specifically, the extraneous cognitive load is imposed by conventional instruction in which the limitation of working memory is rarely taken into account (Kirschner, 2002). The germane cognitive load is defined as the amount of resources devoted to foster the learning process. This type of cognitive load is beneficial, required for the construction and storage of knowledge and information in the long-term memory. Therefore, according to Vogel-Walcutt, et al. (2011), the goal of optimizing cognitive load can be accomplished by: (a) minimizing extraneous cognitive load; (b) maximizing germane cognitive load; and (c) optimizing (increased/decreased as needed) intrinsic cognitive load.

Based on the knowledge of human cognitive architecture and cognitive load theory, Kirschner et al. (2006) present their claim that all of constructivist approaches, included problem-based learning and inquiry learning, are detrimental to learning since they create a huge demand on working memory by pushing learners to search a problem space for problem-relevant information. In other words, when novice learners try to decide what information is important and which information can be considered later or ignored, their lack of the knowledge and experience hinders their ability to distinguish between the two (Vogel-Walcutt, et al., 2011). Furthermore, Kirschner and his colleagues also state that, in fact the knowledge is not stored in the long term memory as the consequences of requiring novice learners to search for problem solutions on their own using a limited working memory.

**Compatibility of Constructivist Approaches with Human Cognitive Architecture**

Some researchers in the domain of cognitive load theory argue that instructional design issues and human cognitive architecture are inseparably intertwined (e.g., Sweller, 2004; Kirschner, 2002). In order to produce an appropriate instructional approach, instructional designers have to take into account the limitation of working memory.

Kirschner and his colleagues (2006) argue that constructivist approaches, in particular problem-based learning, are not likely to be effective because they ignore the findings of cognitive architecture literature that suggest the limits of working memory when dealing with novel information. They further argue that by doing so, constructivist approaches provide learners with partial information that is out of their capability, and thus place a huge burden on working memory when learners are trying to solve problems by searching appropriate information. By allowing learners to construct their own learning experience, the capacity of working memory will be overloaded and learning will be compromised.

However, several constructivist scientists argue that constructivist approaches are in line with human cognitive architecture. In the case of cognitive architecture conceptions, one of proponents of constructivism (Schmidt et al, 2007), states that problem-based learning is compatible with human cognitive architecture because it provides flexible adaptation of guidance either in the level of learner expertise or the complexity of learning task.

Furthermore, Schmidt et al. (2007) present several elements of problem-based learning that can make this approach align with the human cognitive architecture such as,

> (a) students are assembled in *small groups*; (b) these groups receive *training in group collaboration skills* prior to the instruction; (c) their *learning task* is to explain phenomena described in the problem in terms of its underlying principles or mechanism; (d) they do this by initially discussing the problem at hand, *activating* whatever *prior knowledge* is available to each of them; (e) a *tutor* is present to facilitate the learning; (f) s(he) does this by using a *tutor instruction* consisting of relevant information, questions, etc., provided by the problem designer; and (g) *resources* for self-directed study by the students such as books, articles, or other media. (p. 93)

The other researcher, Hmelo-Silver et al. (2007) also suggest that by employing scaffolding extensively in learning process, the cognitive load can be reduced and the learners can learn in more complex domains. These researchers argue that in constructivist approaches such as problem-based learning and inquiry learning, the teachers can provide scaffolding that decrease cognitive load by structuring a task that guide the learners to focus on aspects of the task that are relevant to the learning goals.

In another form of constructivism like discovery learning, Dean and Kuhn (2006) have conducted research about comparing direct instruction and discovery learning with emphasizing in time frame. In this study the researchers found that by allowing students opportunity to develop their strategy in solving problems, they do better than the students who learn by using direct instructional approach. In contrast to the claims made by Kirschner et al. (2006), the students who learned to construct their own strategy lead to significant and lasting gains in strategic understanding than the students who do not (Dean & Kuhn, 2006). This evidence implies that constructivist approaches such as discovery learning can help the learners to acquire knowledge or information and save them in their long-term memory.

**Implication for an Effective Mathematics Teaching: An Example from Realistic Mathematics Education Approach**

Realistic Mathematics Education (RME) approach evolved after 20 years of developmental research at the Freudenthal Institute in Utrecht University, The Netherlands and is thought to have various connections with constructivism. In this sense, I want to point out that RME approach is categorized as constructivism approach. Although both realistic mathematics education approach and constructivism approach share many similarities, there are some differences between them. The constructivism theory is a theory of learning in general, while the realistic mathematics theory is a theory of learning and instruction specifically relating to mathematics. The principles that underlie RME were strongly influenced by Hans Freudenthal's view about mathematics. According to him, mathematics can be best learned by doing it. Education should give the students the "guided" opportunities to be able to reinvent mathematics by doing it themselves. RME attempts to incorporate views on what mathematics is, how students learn mathematics, and how mathematics should be taught.

As in other constructivist approaches, the concept of scaffolding seemed to be present throughout the realist mathematics education theory as well. The teacher in RME classroom must provide scaffolding for students in a way as suggested previously in his lesson plan and must be a facilitator to the students learning. He must never get involved in trying to explain to the pupils. Note that the students in RME classroom are not expected to reinvent everything by themselves. Gravemeijer (1999) emphasizes that the idea is not to motivate students with everyday-life contexts but to look for contexts that are experientially real for students and can be used as starting points for progressive mathematization.

**Conclusion and Suggestion**

I have presented the arguments and evidence from several researches that against the claim of Kirschner et al. (2006) about categorizing constructivist approaches into minimally guidance instruction. The reason is that these constructivist approaches allow for flexible adaptation of guidance and reduce working memory load by giving scaffolding, so that these approaches are compatible with the manner in which our cognitive structures are organized.

Even in this limited review of research on constructivist approaches, it is clear that the claim that "constructivist approaches do not work" is not well supported and in fact, there are supports for the alternative. However, questioning whether constructivist approaches work or not is not really an appropriate question nowadays. Hmelo-Silver, Duncan, and Clark (2007) argue that

> the more important questions to ask are under what conditions do these constructivist approaches work, what kinds of outcomes for which they are effective, what kind of valued practices do they promote, and what kinds of support and scaffolding are needed for different populations and learning goals. (p. 105)

However, I also agree with Kirschner and his colleagues (2006) in the way that our working memory is limited and any instructional approach ignoring this evidence can lead to ineffective way of teaching. I argue that by knowing about working memory load and allowing students to construct their own knowledge under the teacher's scaffolding may lead to be better way of teaching.

I also do not claim that constructivist approaches are better than other instructional guidance such as direct instruction. Instead, I argue that constructivist approaches are compatible with human cognitive architecture as there is flexible adaptation of guidance such as scaffolding in order to avoid working memory load of students. I really agree with Schwartz and his colleagues (1998) that there is a place for both direct instruction and student-directed inquiry. The more challenging question is how to get the balance and sequence right between these two major learning approaches. The claim of Kirschner et al. (2006) that direct instruction is the better instructional approach in which it is appropriate with cognitive load theory is not really well supported. The other study like the one done by Dean and Kuhn (2006) show that direct instruction is neither necessary nor sufficient for robust acquisition or for maintenance over time.

Overall, constructivism, such as realistic mathematics education approach, has had an impressive impact on mathematics education in that it has derived the students into the forefront of doing mathematics. Because of the theory, we realize that careful attention must be paid to how to facilitate students learning in which they are given the opportunities to build their own mathematical knowledge store on the basis of such a learning process.

I wish to conclude this essay with the common wisdom on the nature of mathematics instruction and human learning: "tell me and I will forget; show me and I may remember; involve me and I will understand." The task of education, instead of pouring information into students' heads, is to engage students' minds with powerful and useful concepts.